\documentclass[10pt]{amsart}
\usepackage{hyperref}
\usepackage[utf8]{inputenc}
\usepackage{amsmath,enumitem}
\usepackage{booktabs,url}
\usepackage{amssymb,yhmath,mathrsfs,wasysym}
\usepackage{amsthm}
\usepackage{graphicx}
\usepackage{subfigure}
\usepackage[subfigure]{ccaption}
\usepackage{multicol,multirow,caption}
\usepackage{tikz,pgf}
\usetikzlibrary{arrows,automata,shapes,decorations.pathmorphing,decorations.markings,fadings,patterns}

\title{Statistical Bias in the Distribution of Prime Pairs and Isolated Primes}
\author{Waldemar Puszkarz}
\address{Los Angeles \\ CA \\ USA }
\email{\texttt{psi\_bar@yahoo.com}}
\date{\today}

\begin{document}

\begin{abstract}
Computer experiments reveal that twin primes tend to center on nonsquarefree multiples of 6 more often than on squarefree multiples of 6 compared to what should be expected from the ratio of the number of nonsquarefree multiples of 6 to the number of squarefree multiples of 6 equal $\pi^2/3-1$, or ca 2.290. For multiples of 6 surrounded by twin primes, this ratio is 2.427, a relative difference of ca $6.0\%$ measured against the expected value. A deviation from the expected value of this ratio, ca $1.9\%$, exists also for isolated primes. This shows that the distribution of primes is biased towards nonsquarefree numbers, a phenomenon most likely previously unknown. For twins, this leads to nonsquarefree numbers gaining an excess of $1.2\%$ of the total number of twins. In the case of isolated primes, this excess for nonsquarefree numbers amounts to $0.4\%$ of the total number of such primes. The above numbers are for the first $10^{10}$ primes, with the bias showing a tendency to grow, at least for isolated primes.
\end{abstract}

\maketitle
\thispagestyle{empty}

\section{Surprising effect}
\label{sec:surp}

Except for the first two (2 and 3), all primes are of the form $6k-1$ or $6k+1$. This implies that twin primes (consecutive primes separated by 2) surround multiples of 6 (3 and 5 are the only exception).

All natural numbers are either squarefree or nonsquarefree. Unlike the former, the latter are divisible by a square greater than 1. All primes are squarefree.

Squarefree numbers have natural density of $6/\pi^2$, which gives $1-6/\pi^2$ for natural density of nonsquarefree numbers. Thus, the ratio of relative frequencies at which one expects to find these numbers in a sufficiently large sample of natural numbers is $\pi^2/6-1 = 0.6449$..., which amounts to saying that on average for every 100 squarefree numbers in such a sample there are about 64 nonsquarefree ones.

The ratio in question is different for multiples of 6. Since natural density of squarefree numbers divisible by 6 is $3/\pi^2$, this ratio is $R_{0}=\pi^2/3-1 = 2.2898$... or \textbf{2.290 for practical purposes} (most numbers in this paper are rounded off to the 3rd decimal digit). What this means is that in a sufficiently large sample, on average for every 100 squarefree multiples of 6 there are about 229 nonsquarefree numbers divisible by 6.

If prime pairs are as likely to center on squarefree multiples of 6 as they are on nonsquarefree multiples of 6 (i.e., their distribution is unbiased in this respect), we should expect the same ratio for the multiples of 6 in their centers if calculations are performed on a large enough sample of prime pairs.

\textbf{However, this is not the case.}

This can be observed already in a sample of the first $10^6$ primes and the discrepancy persists (may even be getting slightly stronger) for larger samples. The largest we used consisted of the first $10^{10}$ primes and for it the ratio is $R_{2}=$ \textbf{2.427}.

Defining the relative difference as the absolute value of the ratio of the difference between the experimental value and the theoretical one to the theoretical one, we find out that the relative difference in this case is ca $6.0\%$.

The same calculations applied to isolated primes (primes $p$ such that neither $p-2$ nor $p+2$ is prime) reveal a similar bias:  such primes occur next to nonsquarefree multiples of 6 slightly more often than next to squarefree multiples of 6 compared to what would be expected in the non-biased distribution. In this case, the ratio is $R_{1}=$ \textbf{2.333} (for the sample of $10^{10}$ primes), and the relative difference is ca $1.9\%$.

The bias effect is smaller than for prime pairs, but too large to dismiss it as due to statistical noise.

To be sure this effect has indeed to do with primes and not squarefree numbers in general, we checked if this effect occurs for twin squarefree numbers centered on multiples of 6.

If twin primes are included in such test pairs, our ratio becomes 2.306 for the sample of the first $10^8$ multiples of 6 (and virtually unchanged for the sample of $10^9$), noticeably smaller than 2.427 and very close to the unbiased ratio, 2.290. But if we exclude twin primes from the test pairs, the effect goes away almost completely for the sample of $10^8$ as now we get 2.286, which leads to less than $0.2\%$ in relative difference, a difference small enough to ascribe it to statistical fluctuations. Moreover, for the sample of $10^9$ primes, this ratio is 2.2889, less than $0.05\%$ in relative difference, a very small difference indeed.

For isolated test squarefree numbers next to a multiple of 6 (to the left or right of it), the bias effect pretty much fails to manifest itself already in the sample of the first $10^8$ such multiples as we get 2.2907 (rounded off to the 4th decimal digit), which versus 2.2898 is less than $0.04\%$ in relative difference. If the primes are excluded from these squarefree numbers, we obtain 2.2894, and the relative difference is now ca $0.02\%$, small enough to conclude that also for isolated primes, the bias effect is due to primes; it does not occur for other squarefree numbers. Moreover, for the sample of $10^9$ (primes excluded), the ratio is 2.2897, even closer to the unbiased theoretical value.

\textbf{Hence, to reiterate, the observed effect in both situations is most certainly a property of primes rather than a generic property of squarefree numbers.} To put it more precisely (if not pedantically), if there is any actual contribution to it from non-prime squarefree numbers, it is negligible compared to the contribution from primes.

Moreover, the effect appears pretty stable over several sample ranges with the range size growing by 10 for each data point we collected to determine the effect behavior (see the next section).

\section{Excess functions}
\label{sec:exce}

Let us define them as $\epsilon_{1}=round(1000*(R_{1}-R_{0}))$ for isolated primes and $\epsilon_{2}=round(1000*(R_{2}-R_{0}))$ for prime pairs, where $round(x)$ is tasked with rounding off $x$ to the nearest integer.

The data for $R_{1}$ and $R_{2}$ is obtained numerically while $R_{0}$ can be calculated analytically as done above. We obtained 5 data points for each of these functions (see the data section for more information) and they suggest (albeit quite weakly) that we may be dealing with slowly growing functions, at least in the case of $\epsilon_{1}$. More data is needed to be positive that the trends we suspect these functions may be showing are not due to statistical fluctuations. It may very well be that the best approximation to these functions is a constant. This appears to be the most likely scenario for $\epsilon_{2}$ and it may be so asymptotically for $\epsilon_{1}$.

Below are the values of these functions at arguments that are consecutive powers of 10, starting at $10^6$ (we use the exponent values of $10^{n}$ to index the arguments).

For isolated primes,

$\epsilon_{1}(6)=34$, $\epsilon_{1}(7)=39$, $\epsilon_{1}(8)=41$, $\epsilon_{1}(9)=42$, $\epsilon_{1}(10)=43$.

For prime pairs,

$\epsilon_{2}(6)=136$, $\epsilon_{2}(7)=134$, $\epsilon_{2}(8)=135$, $\epsilon_{2}(9)=136$, $\epsilon_{2}(10)=137$.

These numbers represent the (average) excess of nonsquarefree numbers compared to the non-biased case for every 1000 squarefree numbers. For instance, $\epsilon_{2}(10)$ tells us that there are on average 137 more prime pairs centered on nonsquarefree multiples of 6 per 1000 squarefree multiples of 6 surrounded by primes than one would expect in the unbiased situation for the first $10^{10}$ primes.

\section{Code and data}
\label{sec:code}

The effect discussed was first observed in Mathematica computer experiments performed on the first $10^6$ primes. The data for larger samples was obtained using PARI/GP, an open source software package for number theory.

What follows below is a sample of PARI/GP code used to obtain the data and the data. The data is indexed by exponents of range size, chosen to be powers of 10, from $10^6$ to $10^{10}$. The code for Mathematica can easily be produced from the PARI/GP code.

Our code counts all prime pairs (even though the first of them, $\{3, 5\}$, is not centered on a multiple of 6) and excludes 2 as an isolated prime. While 2 is sometimes treated as an isolated prime, it is actually less isolated from other primes than all odd primes save for 3. With 2 excluded, the number of isolated primes plus twice the number of pairs still add to a range size ($10^6$ through $10^{10}$), for 5 is counted twice as a member of two consecutive pairs that share it. These choices have no impact on our statistical results. We mention them for the sake of clarity.

In what follows, $a$ represents the number of all target objects (primes or test squarefree numbers), while $b$ only the number of such objects next to or centered on squarefree numbers. The ratios discussed above are calculated as $R=(a-b)/b$.

\bigskip
\noindent\textbf{Part A. Prime numbers}

\bigskip
\noindent\textbf{Twin primes}

\bigskip
\begin{verbatim}
a=0; forprime(n=2, prime(10^8), isprime(n+2)&&a++); print1(a)
\end{verbatim}
\begin{verbatim}
\\all twins
\end{verbatim}
\begin{verbatim}
b=0; forprime(n=2, prime(10^8), isprime(n+2)&&issquarefree(n+1)
\end{verbatim}
\begin{verbatim}
&&b++); print1(b) \\twins centered on a squarefree number
\end{verbatim}

\bigskip
\noindent
$a$: 86027 ($6$), 738597 ($7$), 6497407 ($8$), 58047180 ($9$), 524733511 ($10$).

\noindent
$b$: 25113 ($6$), 215732 ($7$), 1897137 ($8$), 16944418 ($9$), 153121114 ($10$).

\bigskip
\noindent\textbf{Isolated primes}

\bigskip
\begin{verbatim}
a=0; forprime(n=3, prime(10^8), !isprime(n+2)&&!isprime(n-2)&&a++);
\end{verbatim}
\begin{verbatim}
print1(a) \\all
\end{verbatim}
\begin{verbatim}
b=0; forprime(n=3, prime(10^8), !isprime(n+2)&&!isprime(n-2)&&
\end{verbatim}
\begin{verbatim}
((n%6==1&&issquarefree(n-1))||(n%6==5&&issquarefree(n+1)))&&b++);
\end{verbatim}
\begin{verbatim}
print1(b) \\next to a squarefree number
\end{verbatim}

\bigskip
\noindent
$a$: 827946 ($6$), 8522806 ($7$), 87005186 ($8$), 883905640 ($9$), 8950532978 ($10$).

\noindent
$b$: 249071 ($6$), 2560208 ($7$), 26123609 ($8$), 265275545 ($9$), 2685404943 ($10$).

\bigskip
\noindent\textbf{Part B. Test squarefree numbers}

\bigskip
\noindent\textbf{Squarefree twins (primes included) centered on a multiple of 6}

\bigskip
\begin{verbatim}
a=0; for(n=1, 10^8, issquarefree(6*n-1)&&issquarefree(6*n+1)&&
\end{verbatim}
\begin{verbatim}
a++); print1(a) \\all
\end{verbatim}
\begin{verbatim}
b=0; for(n=1, 10^8, issquarefree(6*n)&&issquarefree(6*n-1)&&
\end{verbatim}
\begin{verbatim}
issquarefree(6*n+1)&&b++); print1(b) \\centered
\end{verbatim}
\begin{verbatim}
on a squarefree number
\end{verbatim}

\bigskip
\noindent
$a$: 82962973 ($8$), 829630636 ($9$).

\noindent
$b$: 25097397 ($8$), 250974031 ($9$).

\bigskip
\noindent\textbf{Squarefree twins (primes excluded) centered on a multiple of 6}

\bigskip
\begin{verbatim}
a=0; for(n=1, 10^8, issquarefree(6*n-1)&&!isprime(6*n-1)&&
\end{verbatim}
\begin{verbatim}
issquarefree(6*n+1)&&!isprime(6*n+1)&&a++);
\end{verbatim}
\begin{verbatim}
print1(a) \\all
\end{verbatim}
\begin{verbatim}
b=0; for(n=1, 10^8, issquarefree(6*n)&&issquarefree(6*n-1)&&
\end{verbatim}
\begin{verbatim}
!isprime(6*n-1)&&issquarefree(6*n+1)&&!isprime(6*n+1)&&b++);
\end{verbatim}
\begin{verbatim}
print1(b) \\centered on squarefree numbers
\end{verbatim}

\bigskip
\noindent
$a$: 57015536 ($8$), 595982891 ($9$).

\noindent
$b$: 17348734 ($8$), 181210143 ($9$).

\bigskip
\noindent\textbf{Isolated squarefree numbers (primes included) next to a multiple of 6}

\bigskip
\begin{verbatim}
a=0; for(n=1, 10^8, (issquarefree(6*n-1)||issquarefree(6*n+1))&&
\end{verbatim}
\begin{verbatim}
a++); print1(a) \\number of cases a squarefree number is next to
\end{verbatim}
\begin{verbatim}
a multiple of 6
\end{verbatim}
\begin{verbatim}
b=0; for(n=1, 10^8, (issquarefree(6*n))&&(issquarefree(6*n-1)||
\end{verbatim}
\begin{verbatim}
issquarefree(6*n+1))&&b++); print1(b) \\number of cases
\end{verbatim}
\begin{verbatim}
a squarefree number is next to a squarefree multiple of 6
\end{verbatim}

\bigskip
\noindent
$a$: 99415124 ($8$).

\noindent
$b$: 30211331 ($8$).

\bigskip
\noindent\textbf{Isolated squarefree numbers (primes excluded) next to a multiple of 6}

\bigskip
\begin{verbatim}
a=0; for(n=1, 10^8, (issquarefree(6*n-1)&&!isprime(6*n-1))||
\end{verbatim}
\begin{verbatim}
(issquarefree(6*n+1)&&!isprime(6*n+1))&&a++); print1(a) \\number of
\end{verbatim}
\begin{verbatim}
cases a non-prime squarefree number is next to a multiple of 6
\end{verbatim}
\begin{verbatim}
b=0; for(n=1, 10^8, issquarefree(6*n)&&((issquarefree(6*n-1)&&
\end{verbatim}
\begin{verbatim}
!isprime(6*n-1))||(issquarefree(6*n+1)&&!isprime(6*n+1)))&&b++);
\end{verbatim}
\begin{verbatim}
print1(b) \\number of cases a non-prime squarefree number is
\end{verbatim}
\begin{verbatim}
next to a squarefree multiple of 6
\end{verbatim}

\bigskip
\noindent
$a$: 94037859 ($8$), 948253019 ($9$).

\noindent
$b$: 28588317 ($8$), 288245142 ($9$).

\section{Other remarks}
\label{sec:othe}

We mentioned that the effect discussed here can easily be observed even among the first million primes. More surprisingly, though, $R_2$ is consistently larger compared to the expected value even for the first 100, 1000, or 10,000 primes all the way up to one million (and beyond). More than 10,000 first primes were known by the end of the 18th century. What this means is that even Gauss or one of his contemporaries could have noticed it some 200 years ago! Yet, we have found no evidence it was known at all.

The bias that we have so far chosen to measure using the $R$ ratio can also be expressed in other ways. For instance, one can inquire how much the number of twins centered on squarefree numbers, or $b$ (see the data section), differs from its expected value for a given sample size (number $a$ in our data section). Or, we can inquire about the same deviation from the expected value for twins centered on nonsquarefree numbers.

As an example, let us calculate this for the first $10^6$ primes. In this case, we have 86206 prime pairs centered on multiples of 6. If these pairs were distributed in an unbiased way, the way non-prime squarefree twins are, the number of them surrounding squarefree multiples of 6 would be ca $round(86026/(R_{0}+1))=26149$ and not 25113 that we actually get. The deficit we observe, $1036$, represents $3.96\%$ of 26149. This deficit is accompanied by an excess in the number of twins that center on nonsquarefree multiples of 6. Since this excess is $1036$ and the expected number of such twins is $86026-26149=59877$, we get $1.73\%$ for the relative excess measured against the expected value. For isolated primes, similar calculations give $1.03\%$ and $0.45\%$, respectively.

The three measures of bias discussed above are of the same kind. They measure a deviation of some quantity from its expected value relative to this value. The most sensitive of them turns out to be the $R$ ratio, $5.93\%$ and $1.50\%$ for the first $10^6$ primes for pairs and isolated primes, respectively.

One can measure this bias in yet another way suggested to us by Jon E. Schoenfield. This way measures the redistribution of primes due to the bias. For the first $10^6$ primes, $1036$ twins get redistributed compared to the total of $86026$. Thus, as a result of the bias, an excess of $1.20\%$ of the total number of twins is captured by nonsquarefree numbers. For isolated primes, the respective number is $0.31\%$; it grows to $0.39\%$ for the first $10^{10}$ primes.

Since the distribution of non-prime squarefree numbers is not affected by the bias discussed here, one may expect a compensatory effect among the nonsquarefree twins (nonsquarefree numbers surrounding multiples of 6). Such is the case, indeed. The $R$ ratio of these twins is lower than expected. For the first $10^9$ multiples of 6, it is $2.161$, noticeably smaller than $R_0$.

The most common cluster of consecutive squarefree numbers is that of a triple: three squarefree numbers in a row. Doubles and singles are less common, with singles a bit more common than doubles. While this fact may not be widely known, it is not necessarily surprising. But what seems interesting (perhaps even surprising) here is that prime pairs sabotage the formation of squarefree triples by choosing as their middle partners nonsquarefree numbers over squarefree ones more often than is the case among non-prime squarefree twins.

Let us add one more remark that we hope has some pedagogical value. Namely, there is really only one insight invested in this study and rather simple one too: examining the value of $R$. Once the importance of this number is realized, everything else becomes child's play and the whole research can be performed even by a resourceful high school student.

\section{Conclusion}
\label{sec:conc}

The results we reported above are quite basic, were obtained in an elementary fashion, and concern fundamental classes of numbers. It is therefore rather surprising that we found no mention of them in the literature of the subject. This leads us to believe that the statistical bias they describe was most likely unknown.

The study of biases in the distribution of prime numbers has recently been reinvigorated by the work of Lemke Oliver and Soundararajan \cite{Lemke_2016} on the phenomenon related to the one observed by Chebyshev already in 1853 and known as the Chebyshev bias (see \cite{RubinsteinSarnak_1994} and \cite{GranvilleMartin_2006}). However, the effect under consideration here is of different nature than those discussed in the papers cited.

Moreover, and more importantly, the Chebyshev bias is significantly smaller than our bias. Using the data for the first million primes from \cite{Lemke_2016}, we see that the deviation from the non-biased distribution (to be half a million of primes for either of two classes of primes that the Chebyshev effect concerns) is only $170$.

Let us contrast this with the bias presented here. Using the redistribution metric from the previous section, we get a relative bias of $170/1000000$, or less than $0.02\%$, a number over $60$ and $15$ times smaller than the corresponding numbers for twins and isolated primes in our bias.

In statistics, the sample size may matter and comparing findings from samples of considerably different sizes may lead to erroneous conclusions. However, this is not an issue here. Let us note that the size of the sample of twin primes for the first $10^{10}$ primes is 524733511, which makes it comparable to the size of the sample of non-prime squarefree twins for the first $10^9$ multiples of 6, 595982891. Yet, for the latter sample, the ratio $R$ deviates from the expected $R_0$ by less than $0.05\%$, while for the former this deviation is ca $6.0\%$, a two orders of magnitude discrepancy.

It seems that no matter how we look at the data, it is telling us that the bias discussed here is real (and considerably large): primes do occur next to nonsquarefree multiples of 6 more often than next to squarefree multiples of 6 than one would expect from the unbiased distribution, and, in particular, more often than is the case for non-prime squarefree numbers.

One can view this effect as a certain statistical property of primes. It is substantial and simple enough to be of general interest. If only because of this, we believe it deserves further study. The main goal of this paper was to lay the empirical groundwork for this study. More empirical work is still needed, in particular to examine the behavior of excess functions in a wider range of primes to determine these functions approximate analytical form. To do this efficiently, more powerful computing resources are required than we had at our disposal. But what is needed first and foremost is a theoretical model explaining this phenomenon. We hope further research will produce such a model.

\section{Links}
\label{sec:link}

A text file with the PARI/GP code and an Excel spreadsheet with the data and results can be downloaded from \href{http://www.eminimethods.com/Code&Data.zip}{the author's site}.

\bigskip
\noindent\textbf{Acknowledgements.} The author is grateful to the developers of PARI/GP and Wolfram Mathematica, whose software was indispensable to this research, and to Kevin Ford, Krzysztof D. Ma\'slanka, Jon E. Schoenfield, and Marek Wolf for their interest in this work and comments.

\providecommand{\bysame}{\leavevmode\hbox to3em{\hrulefill}\thinspace}
\providecommand{\MR}{\relax\ifhmode\unskip\space\fi MR }
\providecommand{\MRhref}[2]{%
  \href{http://www.ams.org/mathscinet-getitem?mr=#1}{#2}
}
\providecommand{\href}[2]{#2}

\end{document}